\documentclass[12pt,twoside]{article}

\usepackage[english,russian]{babel}
\usepackage{umj_1}           

\usepackage{graphicx, srcltx}   
\usepackage{cite}               

\usepackage{amsmath}            
\usepackage{amsfonts,amssymb}   
\usepackage{srcltx}             

\newcommand{\Leftclt}{Roman A. Veprintsev}
\newcommand{\Rightclt}{On the sharpness of one inequality of different metrics for algebraic polynomials }
\usepackage{authblk}

\DeclareMathOperator*{\esssup}{ess\,sup}

\begin{document}

\hfill \textsc{In memory of my Grandmother}

\begin{Titul}
{\large \bf ON THE SHARPNESS\\ OF ONE INEQUALITY OF DIFFERENT METRICS\\[2mm] FOR ALGEBRAIC POLYNOMIALS }\\[3ex]
{{\bf Roman~A.~Veprintsev} \\[5ex]}
\end{Titul}

\begin{Anot}
{\bf Abstract.} We prove that the previously established inequality of different metrics for algebraic polynomials is sharp in the sense of order.


\end{Anot}


%

\section{Introduction and preliminaries}

In this section, we give some notation used in the article. Here, we also introduce some classes of orthogonal polynomials on $[-1,1]$, including the so-called generalized Gegenbauer polynomials. For a background and more details on the orthogonal polynomials, the reader is referred to \cite{dai_xu_book_approximation_theory_2013,dunkl_xu_book_orthogonal_polynomials_2014,andrews_askey_roy_book_special_functions_1999,szego_book_orthogonal_polynomials_1975}.

Consider the generalized Jacobi weight
\begin{equation*}
\omega_{\alpha,\beta,\gamma}(x)=(1-x)^\alpha(1+x)^\beta|x|^\gamma,\qquad x\in[-1,1],
\end{equation*}
with $\alpha,\,\beta,\,\gamma>-1$. Let $\{p_n^{(\alpha,\beta,\gamma)}\}_{n=0}^\infty$ denote the sequence of orthonormal polynomials with respect to the weight function $\omega_{\alpha,\beta,\gamma}$.

Given $1\leq p\leq\infty$, we denote by $L_p(\omega_{\alpha,\beta,\gamma})$ the space of complex-valued Lebesgue measurable functions $f$ on $[-1,1]$ with finite norm
\begin{equation*}
\begin{array}{lr}
\|f\|_{L_p(\omega_{\alpha,\beta,\gamma})}=\Bigl(\int\nolimits_{-1}^1 |f(x)|^p\,\omega_{\alpha,\beta,\gamma}(x)\,dx\Bigr)^{1/p},&\quad 1\leq p<\infty,\\[1.0em]
\|f\|_{L_\infty(\omega_{\alpha,\beta,\gamma})}=\esssup\limits_{x\in[-1,1]} |f(x)|.&
\end{array}
\end{equation*}

The Jacobi polynomials, denoted by $P_n^{(\alpha,\beta)}(\cdot)$, where $n=0,1,\ldots$, are orthogonal with respect to the Jacobi weight function $w_{\alpha,\beta}(x)=(1-x)^\alpha(1+x)^\beta$ on $[-1,1]$, namely,
\begin{equation*}
\int\nolimits_{-1}^1 P_n^{(\alpha,\beta)}(x)\, P_m^{(\alpha,\beta)}(x)\, w_{\alpha,\beta}(x)\,dt=\begin{cases}\dfrac{2^{\alpha+\beta+1}\Gamma(n+\alpha+1)\Gamma(n+\beta+1)}{(2n+\alpha+\beta+1)\Gamma(n+1)\Gamma(n+\alpha+\beta+1)},&n=m,\\
0,&n\not=m.
\end{cases}
\end{equation*}
Here, as usual, $\Gamma$ is the gamma function. Note that $w_{\alpha,\beta}=\omega_{\alpha,\beta,0}$.

For $\lambda>-\frac{1}{2}$, $\mu\geq0$, and $n=0,\,1,\,\ldots$, the orthonormal generalized Gegenbauer polynomials $\widetilde{C}_n^{(\lambda,\mu)}(\cdot)$ are defined by
\begin{equation}\label{generalized_Gegenbauer_polynomials}
\begin{split}
&\widetilde{C}_{2n}^{(\lambda,\mu)}(x)=\widetilde{a}_{2n}^{\,(\lambda,\mu)}P_n^{(\lambda-1/2,\mu-1/2)}(2x^2-1),\\ &\widetilde{a}_{2n}^{\,(\lambda,\mu)}=\Biggl(\dfrac{(2n+\lambda+\mu)\Gamma(n+1)\Gamma(n+\lambda+\mu)}{\Gamma(n+\lambda+\frac{1}{2})\Gamma(n+\mu+\frac{1}{2})}\Biggr)^{1/2},\\[0.4em]
&\widetilde{C}_{2n+1}^{(\lambda,\mu)}(x)=\widetilde{a}_{2n+1}^{\,(\lambda,\mu)}\,x P_n^{(\lambda-1/2,\mu+1/2)}(2x^2-1),\\ &\widetilde{a}_{2n+1}^{\,(\lambda,\mu)}=\Biggl(\dfrac{(2n+\lambda+\mu+1)\Gamma(n+1)\Gamma(n+\lambda+\mu+1)}{\Gamma(n+\lambda+\frac{1}{2})\Gamma(n+\mu+\frac{3}{2})}\Biggr)^{1/2}.
\end{split}
\end{equation}
They are orthonormal with respect to the weight function
\begin{equation*}
v_{\lambda,\mu}(x)=(1-x^2)^{\lambda-1/2} |x|^{2\mu},\quad x\in[-1,1].
\end{equation*}
Note that $v_{\lambda,\mu}=\omega_{\lambda-1/2,\lambda-1/2,2\mu}$ and
\begin{equation}\label{connection_with_general_polynomials}
\widetilde{C}_{n}^{(\lambda,\mu)}=p_n^{(\lambda-1/2,\lambda-1/2,2\mu)}.
\end{equation}

The generalized Gegenbauer polynomials play an important role in Dunkl harmonic analysis (see, for example, \cite{dunkl_xu_book_orthogonal_polynomials_2014,dai_xu_book_approximation_theory_2013}). So, the study of these polynomials and their applications is very natural.

Define the uniform norm of a continuous function $f$ on $[-1,1]$ by
$
\|f\|_{\infty}=\max\limits_{-1\leq x\leq 1} |f(x)|.
$
The maximum of two real numbers $x$ and $y$ is denoted by $\max(x,y)$.

Throughout the paper we use the following asymptotic notation: $f(n)\lesssim g(n)$, $n\to\infty,$ or equivalently $g(n)\gtrsim f(n)$, $n\to\infty$, means that there exist a positive constant $C$ and a positive integer $n_0$ such that $0\leq f(n)\leq C g(n)$ for all $n\geq n_0$ (asymptotic upper bound); if there exist positive constants $C_1$, $C_2$, and a positive integer $n_0$ such that $0\leq C_1 g(n)\leq f(n)\leq C_2 g(n)$ for all $n\geq n_0$, then we write $f(n)\asymp g(n)$, $n\to\infty$ (asymptotic tight bound).

To simplify the writing, we will omit ``$n\to\infty$'' in the asymptotic notation.

It follows directly from Stirling's asymptotic formula that
\begin{equation}\label{from_Stirling}
\frac{\Gamma(n+\alpha)}{\Gamma(n+\beta)}\asymp n^{\alpha-\beta}
\end{equation}
for arbitrary real numbers $\alpha$ and $\beta$.

It is known \cite[p.~350]{andrews_askey_roy_book_special_functions_1999} that
\begin{equation}\label{uniform_norm_for_Jacobi_polynomials}
\|P_n^{(\alpha,\beta)}\|_{\infty}\asymp
\begin{cases}
n^{\max(\alpha,\beta)}&\text{if }\,\max(\alpha,\beta)\geq-1/2,\\
n^{-1/2}&\text{if }\,\max(\alpha,\beta)<-1/2.
\end{cases}
\end{equation}

The aim of the paper is to prove that the established in \cite{veprintsev_article_inequality_of_different_metrics_2016} inequality of different metrics for algebraic polynomials is sharp in the sense of order.

\section{Main result}

In \cite{veprintsev_article_inequality_of_different_metrics_2016}, we prove the following result.

\begin{teoen}\label{theorem_about_inequality_of_different_metrics}
Let $\alpha\geq\beta\geq-\frac{1}{2}$, $\mu\geq0$, $1\leq p<q\leq\infty$. If $P_n$ is an algebraic polynomial of degree $n$, then
\begin{equation}\label{inequality_of_different_metrics}
\|P_n\|_{L_q(\omega_{\alpha,\beta,\mu})}\lesssim  n^{\max(2(\alpha+1),\mu+1)\left(\frac{1}{p}-\frac{1}{q}\right)} \, \|P_n\|_{L_p(\omega_{\alpha,\beta,\mu})},
\end{equation}
i.e.,
\begin{equation*}
\sup\limits_{P_n} \frac{\|P_n\|_{L_q(\omega_{\alpha,\beta,\mu})}}{\|P_n\|_{L_p(\omega_{\alpha,\beta,\mu})}}\lesssim n^{\max(2(\alpha+1),\mu+1)\left(\frac{1}{p}-\frac{1}{q}\right)},
\end{equation*}
where the supremum is taken over all polynomials $P_n$ of degree $n$.
\end{teoen}

Now we can formulate the main result.

\begin{teoen}\label{main_result}
Let $\alpha\geq\beta\geq-\frac{1}{2}$, $\mu\geq0$, $1\leq p<q\leq\infty$, and $\nu\in\bigl(0,1-\frac{1}{q}\bigr)$. Then
\begin{equation*}
\sup\limits_{P_n} \frac{\|P_n\|_{L_q(\omega_{\alpha,\beta,\mu})}}{\|P_n\|_{L_p(\omega_{\alpha,\beta,\mu})}}\gtrsim \begin{cases}n^{\max(2(\alpha+1),\mu+1)\left(\frac{1}{p}-\frac{1}{q}\right)}&\text{\rm if }\,1<p<q\leq\infty,\\
n^{\max(2(\alpha+1),\mu+1)\left(1-\frac{1}{q}\right)-\nu}&\text{\rm if }\,p=1, \,1<q\leq\infty,
\end{cases}
\end{equation*}
where the supremum is taken over all polynomials $P_n$ of degree $n$.
\end{teoen}

Thus, the inequality~\eqref{inequality_of_different_metrics} of different metrics for algebraic polynomials is precise in order.

\begin{coren}
Let $\alpha\geq\beta\geq-\frac{1}{2}$, $\mu\geq0$, $1<p<q\leq\infty$. Then
\begin{equation*}
\sup\limits_{P_n} \frac{\|P_n\|_{L_q(\omega_{\alpha,\beta,\mu})}}{\|P_n\|_{L_p(\omega_{\alpha,\beta,\mu})}}\asymp n^{\max(2(\alpha+1),\mu+1)\left(\frac{1}{p}-\frac{1}{q}\right)},
\end{equation*}
where the supremum is taken over all polynomials $P_n$ of degree $n$.
\end{coren}

This corollary immediately follows from Theorems~\ref{theorem_about_inequality_of_different_metrics} and \ref{main_result}. 

\section{Some auxiliary results}

In this section, we present some results used in the proof of Theorem~\ref{main_result}.

In \cite{fejzullahu_article_2013} it was proved the following result.

\begin{lemen}
Let $\alpha,\,\beta\geq-1/2$, $\gamma\geq0$, and $\widetilde{\alpha},\,\widetilde{\beta},\,\widetilde{\gamma}>-1$. Let $-1<y_1<0<y_2<1$. Then, for $1\leq p<\infty$,
\begin{equation*}
\int\nolimits_{y_2}^1 (1-x)^{\widetilde{\alpha}} |p_n^{(\alpha,\beta,\gamma)}(x)|^p\,dx\asymp
\begin{cases}
1&\text{\rm if }\,2\widetilde{\alpha}>p\alpha-2+p/2,\\
\ln n&\text{\rm if }\,2\widetilde{\alpha}=p\alpha-2+p/2,\\
n^{p\alpha+p/2-2\widetilde{\alpha}-2}&\text{\rm if }\,2\widetilde{\alpha}<p\alpha-2+p/2,
\end{cases}
\end{equation*}
\begin{equation*}
\int\nolimits_{y_1}^{y_2} |x|^{\widetilde{\gamma}} |p_n^{(\alpha,\beta,\gamma)}(x)|^p\,dx\asymp
\begin{cases}
1&\text{\rm if }\,2\widetilde{\gamma}>p\gamma-2,\\
\ln n&\text{\rm if }\,2\widetilde{\gamma}=p\gamma-2,\\
n^{p\gamma/2-\widetilde{\gamma}-1}&\text{\rm if }\,2\widetilde{\gamma}<p\gamma-2,
\end{cases}
\end{equation*}
\begin{equation*}
\int\nolimits_{-1}^{y_1} (1+x)^{\widetilde{\beta}} |p_n^{(\alpha,\beta,\gamma)}(x)|^p\,dx\asymp
\begin{cases}
1&\text{\rm if }\,2\widetilde{\beta}>p\beta-2+p/2,\\
\ln n&\text{\rm if }\,2\widetilde{\beta}=p\beta-2+p/2,\\
n^{p\beta+p/2-2\widetilde{\beta}-2}&\text{\rm if }\,2\widetilde{\beta}<p\beta-2+p/2.
\end{cases}
\end{equation*}
\end{lemen}

Using this lemma and \eqref{connection_with_general_polynomials}, one can deduce the following two corollaries.

\begin{coren}\label{first_main_corollary}
Let $\alpha\geq\beta\geq-\frac{1}{2}$, $\mu\geq0$, and $1<p<\infty$. Then
\begin{equation*}
\|p_n^{(2\alpha+\frac{3}{2},2\beta+\frac{3}{2},2\mu+2)}\|_{L_p(\omega_{\alpha,\beta,\mu})}\asymp n^{\max(2(\alpha+1),\mu+1)\left(1-\frac{1}{p}\right)}.
\end{equation*}
In particular,
\begin{equation*}
\|\widetilde{C}_n^{(2(\alpha+1),\mu+1)}\|_{L_p(\omega_{\alpha,\beta,\mu})}\asymp
n^{\max(2(\alpha+1),\mu+1)\left(1-\frac{1}{p}\right)}.
\end{equation*}
\end{coren}

\begin{coren}\label{second_main_corollary}
Let $\alpha\geq\beta\geq-\frac{1}{2}$, $\mu\geq0$, $\nu\in\bigl(0,1-\frac{1}{q}\bigr)$, and $1<q<\infty$. Then
\begin{equation*}
\begin{split}
&\|p_n^{(2\alpha+\frac{3}{2}-\nu,2\beta+\frac{3}{2}-\nu,2\mu+2-2\nu)}\|_{L_1(\omega_{\alpha,\beta,\mu})}\asymp1,\\
&\|p_n^{(2\alpha+\frac{3}{2}-\nu,2\beta+\frac{3}{2}-\nu,2\mu+2-2\nu)}\|_{L_q(\omega_{\alpha,\beta,\mu})}\asymp n^{\max(2(\alpha+1),\mu+1)\left(1-\frac{1}{q}\right)-\nu}.
\end{split}
\end{equation*}
In particular,
\begin{equation*}
\|\widetilde{C}_n^{(2(\alpha+1)-\nu,\mu+1-\nu)}\|_{L_1(\omega_{\alpha,\beta,\mu})}\asymp1.
\end{equation*}
\end{coren}

For $\lambda>-\frac{1}{2}$, $\mu>0$, we have \cite[Theorem~1]{veprintsev_article_maximum_2016}
\begin{equation}\label{max_value_of_generalized_Gegenbauer_polynomials}
\|\widetilde{C}_n^{(\lambda,\mu)}\|_{\infty}\asymp n^{\max(\lambda,\mu)}.
\end{equation}
The above formula is sufficient for our purpose, but let us show that this is valid for $\mu=0$.

\begin{prpen}
Let $\lambda>-\frac{1}{2}$. Then
\begin{equation*}
\|\widetilde{C}_n^{(\lambda,0)}\|_{\infty}\asymp n^{\max(\lambda,0)}.
\end{equation*}
\end{prpen}

\proofen From \eqref{generalized_Gegenbauer_polynomials} and \cite[Theorem~4.1]{szego_book_orthogonal_polynomials_1975}, we have
\begin{equation*}
\begin{split}
&\widetilde{C}_{2n}^{(\lambda,0)}(x)=\widetilde{a}_{2n}^{\,(\lambda,0)} \, P_n^{(\lambda-1/2,-1/2)}(2x^2-1)=\\
&\hspace{3.8em}=\widetilde{a}_{2n}^{\,(\lambda,0)} \,\frac{\Gamma(n+\lambda+\frac{1}{2})\, \Gamma(2n+1)}{\Gamma(2n+\lambda+\frac{1}{2})\,\Gamma(n+1)} \, P_{2n}^{(\lambda-1/2,\lambda-1/2)}(x),\\
&\widetilde{C}_{2n+1}^{(\lambda,0)}(x)=\widetilde{a}_{2n+1}^{\,(\lambda,0)} \, x P_n^{(\lambda-1/2,1/2)}(2x^2-1)=\\
&\hspace{3.8em}=\widetilde{a}_{2n+1}^{\,(\lambda,0)} \,\frac{\Gamma(n+\lambda+\frac{1}{2})\, \Gamma(2n+2)}{\Gamma(2n+\lambda+\frac{3}{2})\,\Gamma(n+1)} \, P_{2n+1}^{(\lambda-1/2,\lambda-1/2)}(x).
\end{split}
\end{equation*}
By \eqref{generalized_Gegenbauer_polynomials} and \eqref{from_Stirling},
\begin{equation*}\label{asymptotic_for_coefficients_for_orthonormal_generalized_Gegenbauer_polynomials}
\widetilde{a}_{2n}^{\,(\lambda,0)}\asymp n^{1/2},\qquad \widetilde{a}_{2n+1}^{\,(\lambda,0)}\asymp n^{1/2},
\end{equation*}
\begin{equation*}
\frac{\Gamma(n+\lambda+\frac{1}{2})\, \Gamma(2n+1)}{\Gamma(2n+\lambda+\frac{1}{2})\,\Gamma(n+1)}\asymp 1,\qquad \frac{\Gamma(n+\lambda+\frac{1}{2})\, \Gamma(2n+2)}{\Gamma(2n+\lambda+\frac{3}{2})\,\Gamma(n+1)}\asymp 1.
\end{equation*}
Hence, considering \eqref{uniform_norm_for_Jacobi_polynomials}, we get the desired asymptotic behavior
\begin{equation*}
\begin{split}
\|\widetilde{C}_n^{(\lambda,0)}\|_\infty&\asymp n^{1/2} \|P_n^{(\lambda-1/2,\lambda-1/2)}\|_\infty\asymp\\
&\asymp n^{1/2} \cdot
\begin{cases}
n^{\lambda-1/2}&\text{if }\,\lambda\geq0,\\
n^{-1/2}&\text{if }\,-1/2<\lambda<0.
\end{cases}
\end{split}
\end{equation*}
\hfill$\square$

\section{Proof of Theorem~\ref{main_result}}

Let us consider the following cases:
\begin{itemize}
\item[I)] $1<p<q<\infty$ and $1<p<\infty$, $q=\infty$;
\item[II)] $p=1$, $1<q<\infty$ and $p=1$, $q=\infty$.
\end{itemize}

Case I). By Corollary \ref{first_main_corollary} and \eqref{max_value_of_generalized_Gegenbauer_polynomials}, we get the desired behavior
\begin{equation*}
\begin{split}
\sup\limits_{P_n} \frac{\|P_n\|_{L_q(\omega_{\alpha,\beta,\mu})}}{\|P_n\|_{L_p(\omega_{\alpha,\beta,\mu})}}&\geq \frac{\|p_n^{(2\alpha+\frac{3}{2},2\beta+\frac{3}{2},2\mu+2)}\|_{L_q(\omega_{\alpha,\beta,\mu})}}{\|p_n^{(2\alpha+\frac{3}{2},2\beta+\frac{3}{2},2\mu+2)}\|_{L_p(\omega_{\alpha,\beta,\mu})}}\gtrsim\\
&\gtrsim \frac{n^{\max(2(\alpha+1),\mu+1)\left(1-\frac{1}{q}\right)}}{n^{\max(2(\alpha+1),\mu+1)\left(1-\frac{1}{p}\right)}}=n^{\max(2(\alpha+1),\mu+1)\left(\frac{1}{p}-\frac{1}{q}\right)}
\end{split}
\end{equation*}
and
\begin{equation*}
\begin{split}
\sup\limits_{P_n} \frac{\|P_n\|_{\infty}}{\|P_n\|_{L_p(\omega_{\alpha,\beta,\mu})}}&\geq  \frac{\|\widetilde{C}_n^{(2(\alpha+1),\mu+1)}\|_\infty}{\|\widetilde{C}_n^{(2(\alpha+1),\mu+1)}\|_{L_p(\omega_{\alpha,\beta,\mu})}}\gtrsim\\
&\gtrsim\frac{n^{\max(2(\alpha+1),\mu+1)}}{n^{\max(2(\alpha+1),\mu+1)\left(1-\frac{1}{p}\right)}}=n^{\max(2(\alpha+1),\mu+1)\frac{1}{p}}.
\end{split}
\end{equation*}

Case II). Using Corollary \ref{second_main_corollary} and \eqref{max_value_of_generalized_Gegenbauer_polynomials}, we obtain the following asymptotic estimates
\begin{equation*}
\begin{split}
\sup\limits_{P_n} \frac{\|P_n\|_{L_q(\omega_{\alpha,\beta,\mu})}}{\|P_n\|_{L_1(\omega_{\alpha,\beta,\mu})}}&\geq \frac{\|p_n^{(2\alpha+\frac{3}{2}-\nu,2\beta+\frac{3}{2}-\nu,2\mu+2-2\nu)}\|_{L_q(\omega_{\alpha,\beta,\mu})}}{\|p_n^{(2\alpha+\frac{3}{2}-\nu,2\beta+\frac{3}{2}-\nu,2\mu+2-2\nu)}\|_{L_1(\omega_{\alpha,\beta,\mu})}}\gtrsim\\
&\gtrsim n^{\max(2(\alpha+1),\mu+1)\left(1-\frac{1}{q}\right)-\nu}
\end{split}
\end{equation*}
and
\begin{equation*}
\begin{split}
\sup\limits_{P_n} \frac{\|P_n\|_{\infty}}{\|P_n\|_{L_1(\omega_{\alpha,\beta,\mu})}}&\geq\frac{\|\widetilde{C}_n^{(2(\alpha+1)-\nu,\mu+1-\nu)}\|_{\infty}}{\|\widetilde{C}_n^{(2(\alpha+1)-\nu,\mu+1-\nu)}\|_{L_1(\omega_{\alpha,\beta,\mu})}}\gtrsim\\
&\gtrsim n^{\max(2(\alpha+1),\mu+1)-\nu}.
\end{split}
\end{equation*}

\begin{Biblioen}

\bibitem{andrews_askey_roy_book_special_functions_1999}G.\,E.~Andrews, R.~Askey, and R.~Roy, \textit{Special Functions}, Encyclopedia of Mathematics and its Applications \textbf{71}, Cambridge University Press, Cambridge, 1999.

\bibitem{dai_xu_book_approximation_theory_2013}F.~Dai and Y.~Xu, \textit{Approximation theory and harmonic analysis on spheres and balls}, Springer Monographs in Mathematics, Springer, 2013.

\bibitem{dunkl_xu_book_orthogonal_polynomials_2014}C.\,F.~Dunkl and Y.~Xu, \textit{Orthogonal polynomials of several variables}, 2nd ed., Encyclopedia of Mathematics and its Applications \textbf{155}, Cambridge University Press, Cambridge, 2014.

\bibitem{fejzullahu_article_2013}B.\,Xh.~Fejzullahu, On orthogonal expansions with respect to the generalized Jacobi weight, \textit{Results. Math.} \textbf{63}~(2013), 1177--1193.

\bibitem{szego_book_orthogonal_polynomials_1975}G.~Szeg\"{o}, \textit{Orthogonal polynomials}, 4th ed., American Mathematical Society Colloquium Publications \textbf{23}, American Mathematical Society, Providence, Rhode Island, 1975.

\bibitem{veprintsev_article_maximum_2016}R.\,A.~Veprintsev, On the asymptotic behavior of the maximum absolute value of generalized Gegenbauer polynomials, arXiv preprint 1602.01023 (2016).

\bibitem{veprintsev_article_inequality_of_different_metrics_2016}R.\,A.~Veprintsev, On an inequality of different metrics for algebraic polynomials, arXiv preprint 1606.06149 (2016).
\end{Biblioen}

\noindent \textsc{Independent researcher, Uzlovaya, Russia}

\noindent \textit{E-mail address}: \textbf{veprintsevroma@gmail.com}

\end{document}